\documentclass[10pt]{article}
\begin{document}
\title{ {\bf An Extremal Problem On Potentially $K_{m}-C_{4}$-graphic
Sequences}
\thanks{  Project Supported by NNSF of China(10271105), NSF of Fujian(Z0511034),
Science and Technology Project of Fujian, Fujian Provincial
Training Foundation for "Bai-Quan-Wan Talents Engineering" ,
Project of Fujian Education Department and Project of Zhangzhou
Teachers College.}}
\author{{Chunhui Lai}\\
{\small Department of Mathematics}\\{\small Zhangzhou Teachers
College, Zhangzhou} \\{\small Fujian 363000,
 P. R. of CHINA.}\\{\small e-mail: zjlaichu@public.zzptt.fj.cn}}
\date{}
\maketitle
\begin{center}
\begin{minipage}{120mm}
\vskip 0.1in
\begin{center}{\bf Abstract}\end{center}
 {A sequence $S$ is potentially $K_{m}-C_{4}$-graphical if it has
a realization containing a $K_{m}-C_{4}$ as a subgraph. Let
$\sigma(K_{m}-C_{4}, n)$ denote the smallest degree sum such that
every $n$-term graphical sequence $S$ with $\sigma(S)\geq
\sigma(K_{m}-C_{4}, n)$ is potentially $K_{m}-C_{4}$-graphical. In
this paper, we prove that $\sigma (K_{m}-C_{4}, n)\geq
(2m-6)n-(m-3)(m-2)+2,$ for $n \geq m \geq 4.$  We conjecture that
equality holds for $n \geq m \geq 4.$ We prove that this
conjecture is true for $m=5$.}\par
\par
 {\bf Key words:} graph; degree sequence; potentially $K_{m}-C_{4}$-graphic
sequence\par
  {\bf AMS Subject Classifications:} 05C07, 05C35\par
\end{minipage}
\end{center}
 \par
 \section{Introduction}
\par

  If $S=(d_1,d_2,...,d_n)$ is a sequence of
non-negative integers, then it is called  graphical if there is a
simple graph $G$ of order $n$, whose degree sequence ($d(v_1 ),$
$d(v_2 ),$ $...,$ $d(v_n )$) is precisely $S$. If $G$ is such a
graph then $G$ is said to realize $S$ or be a realization of $S$.
A graphical sequence $S$ is potentially $H$-graphical if there is
a realization of $S$ containing $H$ as a subgraph, while $S$ is
forcibly $H$-graphical if every realization of $S$ contains $H$ as
a subgraph. Let $\sigma(S)=d(v_1 )+d(v_2 )+... +d(v_n ),$ and
$[x]$ denote the largest integer less than or equal to $x$. We
denote  $G+H$ as the graph with $V(G+H)=V(G)\bigcup V(H)$ and
$E(G+H)=E(G)\bigcup E(H)\bigcup \{xy: x\in V(G) , y \in V(H) \}. $
Let $K_k$, and $C_k$ denote a complete graph on $k$ vertices, and
a cycle on $k$ vertices, respectively. Let $K_{m}-C_{4}$ be the
graph obtained from $K_{m}$ by removing four edges of a 4 cycle
$C_{4}$.
\par

Given a graph $H$, what is the maximum number of edges of a graph
with $n$ vertices not containing $H$ as a subgraph? This number is
denoted $ex(n,H)$, and is known as the Tur\'{a}n number. This
problem was proposed for $H = C_4$ by Erd\"os [2] in 1938 and in
general by Tur\'{a}n [10]. In terms of graphic sequences, the
number $2ex(n,H)+2$ is the minimum even integer $l$ such that
every $n$-term graphical sequence $S$ with $\sigma (S)\geq l $ is
forcibly $H$-graphical. Here we consider the following variant:
determine the minimum even integer $l$ such that every $n$-term
graphical sequence $S$ with $\sigma(S)\ge l$ is potentially
$H$-graphical. We denote this minimum $l$ by $\sigma(H, n)$.
Erd\"os,\ Jacobson and Lehel [3] showed that $\sigma(K_k, n)\ge
(k-2)(2n-k+1)+2$ and conjectured that equality holds. They proved
that if $S$ does not contain zero terms, this conjecture is true
for $k=3,\ n\ge 6$. The conjecture is confirmed in [4],[6],[7],[8]
and [9].
 \par
 Gould,\ Jacobson and
Lehel [4] also proved that  $\sigma(pK_2, n)=(p-1)(2n-2)+2$ for
$p\ge 2$; $\sigma(C_4, n)=2[{{3n-1}\over 2}]$ for $n\ge 4$.  Lai
[5] proved that  $\sigma (K_4-e, n)=2[{{3n-1}\over 2}]$ for $n\ge
7$.\  In this paper, we prove that $\sigma (K_{m}-C_{4}, n)\geq
(2m-6)n-(m-3)(m-2)+2,$ for $n \geq m \geq 4.$  We conjectured that
equality holds for $n \geq m \geq 4.$ We prove that this
conjecture is true for $m=5$.\par

\section{ Main results.} \par
{\bf  Theorem 1.} $\sigma (K_{m}-C_{4}, n)\geq
(2m-6)n-(m-3)(m-2)+2,$ for $n \geq m \geq 4.$
\par
{\bf Proof.}   Let
$$H=K_{m-3}+  \overline{K_{n-m+3}}$$
Then $H$ is a uniquely realization of $((n-1)^{m-3},
(m-3)^{n-m+3})$ and $H$ clearly does not contain $K_{m}-C_{4}.$
Thus
$$\sigma (K_{m}-C_{4}, n) \geq (m-3)(n-1) + (m-3)(n-m+3)+ 2
= (2m-6)n-(m-3)(m-2)+2.$$
\par
{\bf  Theorem 2.} For  $n\geq 5$,
    $\sigma (K_{5}-C_{4}, n)=4n-4.$

\par
{\bf Proof.} By theorem 1, for $n\geq5,$ $\sigma (K_{5}-C_{4}, n)
\geq 4n-4.$ We need to show that if $S$ is an $n$-term graphical
sequence with $\sigma(S)\geq 4n-4$, then there is a realization of
$S$ containing a $K_{5}-C_{4}.$
 Let $d_{1} \geq d_{2} \geq \cdots \geq d_{n}$, and let $G$ be a realization of $S.$
 \par
   Case: $n=5$, if a graph has size
$q \geq 8$, then clearly it contains a $K_{5}-C_{4}$, so that
$\sigma(K_{5}-C_{4} ,5)\leq 4n-4$.
\par
  Case:  $n=6$. If $\sigma(S)=20$, we first consider $d_{6} \leq2$. Let
  $S^{'}$ be the degree sequence of $G-v_{6}$, so
  $\sigma(S^{'})\geq20-2\times2=16$. Then, by induction $S^{'}$ has a realization
  containing a $K_{5}-C_{4}$. Therefore $S$ has a realization
   containing a $K_{5}-C_{4}$. Now we consider $d_{6}\geq 3$. It is
  easy to see that $S=(5^{1},3^{5})$
  or $S=(4^{2},3^{4})$. Obviously, each is potentially
  $K_{5}-C_{4}$-graphic.
   Next, if $\sigma(S)=22$ then it must be that $d_{6}\leq3$. Let
  $ S^{'}$ be the degree sequence of $G-v_{6}$, so
   $\sigma(S^{'})\geq22-3\times2=16$. Then $S^{'}$ has a realization
   containing a $K_{5}-C_{4}$. Therefore $S$ has a realization
   containing a $K_{5}-C_{4}$.
   Finally, suppose that $\sigma(S)\geq24$. We first consider
   $d_{6}\leq4$. Let
 $ S^{'}$ be the degree sequence of $G-v_{6}$, so
  $\sigma(S^{'})\geq24-2\times4=16$. Then $S^{'}$ has a realization
  containing a $K_{5}-C_{4}$. Therefore $S$ has a realization
   containing a $K_{5}-C_{4}$. Now we consider $d_{6}\geq5$. It is
   easy to see that $S=(5^{6})$. Obviously, $(5^{6})$ is potentially $K_{5}-C_{4}$-graphic.

  \par
   Case:  $n=7$. First we assume that $\sigma(S)=24$. Suppose $d_{7}\leq2$ and let $S^{'}$ be the degree sequence of
    $G-v_{7}$, so $\sigma(S^{'})\geq24-2\times2=20$. Then $S^{'}$ has a realization
  containing a $K_{5}-C_{4}.$  Therefore $S$ has a realization
   containing a $K_{5}-C_{4}.$   Now we
   assume that $d_{7}\geq3$. It is easy to see that $S$ is one of
   $(6^{1},3^{6})$,
    $(5^{1},4^{1},3^{5})$, or $(4^{3},3^{4}).$
   Obviously, all of them are  potentially $K_{5}-C_{4}$-graphic.
Next,  if $\sigma(S)=26$, It is
   easy to see that $d_{7}\leq3$. Let $S^{'}$ be the
degree sequence of $G-v_{7}$, so
$\sigma(S^{'})\geq26-3\times2=20$. Then $S^{'}$ has a realization
  containing a $K_{5}-C_{4}$. Therefore $S$ has a realization
   containing a $K_{5}-C_{4}.$
  Finally, suppose that $\sigma\geq28$. If $d_{7}\leq4$. Let $S^{'}$ be the degree sequence
  of $G-v_{7}$, so $\sigma(S^{'})\geq28-2\times4=20$. Then $S^{'}$ has a realization
  containing a $K_{5}-C_{4}$. Therefore $S$ has a realization
   containing a $K_{5}-C_{4}$.
      Now we consider $d_{7}\geq5$. It is easy to see that
   $\sigma(S)>5\times7=35$.
    Clearly,  $d_{7}\leq6$. Let $S^{'}$ be the
degree sequence of $G-v_{7}$, so $\sigma(S^{'})>35-6\times2=23$.
Then, by induction $S^{'}$ has a realization
  containing a $K_{5}-C_{4}$. Therefore $S$ has a realization
   containing a $K_{5}-C_{4}.$
   \par
   We proceed by induction on $n$. Take $n \geq 8$ and make the
   inductive assumption that for $7 \leq t < n$, whenever $S_{1}$ is
   a $t$-term graphical sequence such that
          $$\sigma(S_{1})\geq4t-4$$
   then $S_{1}$ has a realization containing a $K_{5}-C_{4}.$
      Let $S$ be an $n$-term graphical sequence with $\sigma(S)\geq4n-4$.
      If $d_{n}\leq2$, let $S^{'}$ be the degree sequence of $G-v_{n}$. Then
   $\sigma(S^{'})\geq4n-4-2\times2=4(n-1)-4$. By induction, $S^{'}$ has a realization
   containing a $K_{5}-C_{4}$. Therefore $S$ has a realization containing a
   $K_{5}-C_{4}$. Hence, we may assume that $d_{n}\geq3$. By
   Proposition 2 and Theorem 4 of [4] (or Theorem  3.3 of [6] ) $S$ has a
   realization containing a $K_{4}$. By Lemma 1 of [4] ,there is a
   realization $G$ of $S$ with $v_{1},v_{2},v_{3},v_{4}$, the
   four vertices of highest degree containing a $K_{4}$.
     If $d(v_{2})=3$, then $4n-4 \leq \sigma(S) \leq n-1+3(n-1)=4n-4$.
     Hence, $S= ((n-1)^{1}, 3^{n-1})$. Obviously, $((n-1)^{1}, 3^{n-1})$ is potentially
  $K_{5}-C_{4}$-graphic. Therefore, we may assume that $d(v_{2}) \geq 4$.
      Let $v_{1}$ be adjacent to $v_{2},v_{3},v_{4},y_{1}$. If $y_{1}$ is adjacent
   to one of $v_{2},v_{3},v_{4}$, then $G$ contains a $K_{5}-C_{4}$. Hence, we may assume
   that $y_{1}$ is not adjacent to $v_{2},v_{3},v_{4}$.
      Let $v_{2}$ be adjacent to $v_{1},v_{3},v_{4},y_{2}$. If $y_{2}$ is adjacent
   to one of $v_{1},v_{3},v_{4}$, then $G$ contains a $K_{5}-C_{4}$. Hence, we may assume
   that $y_{2}$ is not adjacent to $v_{1},v_{3},v_{4}$.
       Since $d(y_{1}) \geq d_{n} \geq 3$, there is a new vertex $y_{3}$, such that
   $y_{1}y_{3} \in E(G)$.
   \par

       If $y_{3}v_{1} \in E(G)$, then $G$ contains a $K_{5}-C_{4}$. Hence, we may assume
   that $y_{3}v_{1}  \notin  E(G)$. Then the edge interchange that removes the edges
   $y_{1}y_{3},v_{1}v_{4}$ and $v_{2}y_{2}$ and inserts the non-edges
   $y_{1}v_{2},y_{3}v_{1}$ and $y_{2}v_{4}$ produces a realization
   $G^{'}$ of $S$ containing a $K_{5}-C_{4}$ where $K_{5}-C_{4}$ with five vertices $v_{1},v_{2},v_{3},v_{4},y_{1}$.
       \par

       This finishes the inductive step, and thus Theorem 2 is established.
      \par
       We make the following conjecture:
       \par
     {\bf Conjecture.} $$\sigma (K_{m}-C_{4}, n)=
(2m-6)n-(m-3)(m-2)+2$$ for $n \geq m \geq 4.$
   \par
       This conjecture is true for $m=4$, by Theorem 5  of [4], and
       for $m=5$, by the above Theorem.
       \par

        \section*{Acknowledgment}
  The author thanks the referees for many helpful comments.

\end{document}